\documentclass[11pt,twoside,english]{myamsart}
\advance\oddsidemargin by -1.0cm
\advance\evensidemargin by -1.0cm
\advance\topmargin by -1.0cm 
\usepackage{amssymb}
\usepackage{babel}
\usepackage{amstext}
\usepackage{amscd}   
\usepackage{epsfig}  
\usepackage{rotating}

\theoremstyle{plain}
\newtheorem{cor}{Corollary}[section]
\newtheorem{lem}{Lemma}[section]
\newtheorem{thm}{Theorem}[section]            
\newtheorem{prop}{Proposition}[section]
\theoremstyle{definition}
\newtheorem{exa}{Example}[section]
\newtheorem{NB}{Remark}[section]

\renewcommand{\labelenumi}{$(\theenumi)$}
\newcommand{\bdm}{\begin{displaymath}}
\newcommand{\edm}{\end{displaymath}}
\newcommand{\be}{\begin{equation}}
\newcommand{\ee}{\end{equation}}
\newcommand{\ba}[1]{\begin{array}{#1}}
\newcommand{\ea}{\end{array}}

\newcommand{\btab}{\begin{tabular}}
\newcommand{\etab}{\end{tabular}}

\newcommand{\R}{\ensuremath{\mathbb{R}}}

\newcommand{\G}{\ensuremath{\mathrm{G}}}

\newcommand{\SU}{\ensuremath{\mathrm{SU}}}

\newcommand{\SO}{\ensuremath{\mathrm{SO}}}
\newcommand{\Spin}{\ensuremath{\mathrm{Spin}}}

\begin{document}
\def\haken{\mathbin{\hbox to 6pt{%
                 \vrule height0.4pt width5pt depth0pt
                 \kern-.4pt
                 \vrule height6pt width0.4pt depth0pt\hss}}}
    \let \hook\intprod
\setcounter{equation}{0}
%
%
\thispagestyle{empty}
%
\date{\today}
\title[The second Dirac eigenvalue of a nearly parallel 
$\G_2$-manifold]
{ The second Dirac eigenvalue of a nearly parallel 
$\G_2$-manifold }
%
%
%
\author{Thomas Friedrich}
\address{\hspace{-5mm}  
Thomas Friedrich\newline
Institut f\"ur Mathematik \newline
Humboldt-Universit\"at zu Berlin\newline
Sitz: WBC Adlershof\newline
D-10099 Berlin, Germany\newline
{\normalfont\ttfamily friedric@mathematik.hu-berlin.de}}
%
\subjclass[2000]{Primary 53 C 25; Secondary 81 T 30}
\begin{abstract}
We investigate the second Dirac eigenvalue on Riemannian manifolds
admitting a Killing spinor. In small dimensions the whole Dirac spectrum
depends on special eigenvalues on functions and $1$-forms. We
compute  and discuss the formulas in dimension $n = 7$.
\end{abstract}
\maketitle
\pagestyle{headings}
%
%
\section{Introduction}\noindent
%
Let $(M^n,g)$ be a compact Riemannian spin manifold admitting a Killing spinor
$\psi$,
\bdm
\nabla_X \psi = \ a \cdot X \cdot \psi \, , \quad a \in \R^1 \ . 
\edm
The $M^n$ is an Einstein space with non-negative scalar curvature $R \geq 0$ and
\bdm
n^2 \, a^2 \ = \ \mu_1(D^2) \, = \ \frac{n}{4(n-1)} R
\edm
is the smallest eigenvalues of the square of the Riemannian Dirac operator $D$, see
\cite{Fri1}. The question whether or not one can estimate the next
eigenvalue $\mu_2(D^2)$ has not yet been investigated for Dirac operators.
Remark that in case of the Laplacian acting on functions of an Einstein space
$M^n \neq S^n$, there are lower estimates  for small eigenvalues
depending on the minimum of the
sectional curvature, see \cite{Simon}, \cite{Tanno}.\\ 

\noindent
The Killing
spinor of $M^n$
allows us to construct out of it other spinor fields. We can use, for example,
test spinors $\psi^* = f \cdot \psi + c \cdot df \cdot \psi$, 
where $f$ is an eigenfunction of the
Laplace operator. In this way we obtain an  {\it upper bound} depending on
the first positive eigenvalue  $\lambda^0_1$ of the Laplace operator
on functions. A more difficult question is to find {\it lower bounds} for
$\mu_2(D^2)$ on Riemannian spin manifolds with Killing spinors. We study the
problem  comparing  the Dirac spectrum
with the Laplace spectrum on functions and $1$-forms. In more details we will
discuss the question in dimension 7 and remark that this method
should work in dimensions $n=5$ ($5$-dimensional Sasaki-Einstein manifolds) and
$n= 6$ (nearly K\"ahler manifolds), too. In dimension $n =3$ a similar
question has been discussed recently by E.C. Kim, see \cite{Kim}.

\section{Estimates for $\mu_2(D^2)$ in arbitrary dimension }\noindent
%
Consider an eigenspinor $\psi^*$ of the Dirac operator
\bdm
D(\psi^*) \ = \ m \cdot \psi^* \ , \quad \int_{M^n} \langle \psi \, , \psi^* 
\rangle \ =
\ 0
\edm
being $L^2$-orthogonal to the Killing spinor $\psi$. Moreover, suppose that the
function $f := \langle \psi \, , \, \psi^* \rangle$ is not identically zero. A direct
computation yield the formula
\bdm
\int_{M^n} \Delta(f) \cdot f \ = \ \Big\{ m^2 \, + \, 2 a m \, + \, a^2
(2n \, - \, n^2) \Big\} \int_{M^n} f^2 \ .
\edm 
Thus we obtain the estimate
\bdm
\lambda^0_1 \ \leq \ m^2 \, + \, 2 a m \, + \, a^2
(2n \, - \, n^2) \ , 
\edm
or,  respectively, 
\bdm
\sqrt{\lambda^0_1 \, + \, a^2 (1 - n)^2} \ - \ |a| \ \leq \ |m| \ .
\edm
Conversely, if $f$ is a non-trivial eigenfunction, $\Delta(f) = \lambda^0_1
f$, then the two numbers
\bdm
m \ := \ - \, a \, \pm \, \sqrt{\lambda^0_1 \, + \, a^2 (1 - n)^2} 
\edm
are eigenvalues of the Dirac operator. The corresponding 
eigenspinor is given by the formula
\bdm
\psi^* \ := \ f \cdot \psi \, + \, \frac{1}{m + 2a - na} df \cdot \psi \ .
\edm
\begin{NB}
If $M^n \neq S^n$ is not isometric to the sphere, then the Lichnerowicz-Obata
theorem estimates $\lambda^0_1$,
\bdm
\lambda^0_1 \ > \ \frac{R}{n-1} \ = \ 4 a^2 n.
\edm
Therefore our lower bound 
\bdm
\sqrt{\lambda^0_1 \, + \, a^2 (1 - n)^2} \ - \ |a| \ > \ \sqrt{\mu_1(D^2)} \ = \
|a| \, n 
\edm
is greater then the smallest eigenvalue of the Dirac operator.
\end{NB}
\noindent
Let us summarize the result.
\begin{thm}
Let  $M^n \neq S^n$ be a compact Riemannian spin manifold with a Killing spinor
$\psi$, $\nabla_X \psi = a \cdot X \cdot \psi$. Then the first
eigenvalue of the square $D^2$ of the Dirac operator equals  $\mu_1(D^2) = a^2
n^2$. The numbers
\bdm
 \Big(\pm \, \sqrt{\lambda^0_1 \, + \, a^2 (1 - n)^2} \ - \ |a|\Big)^2 
\edm
are  eigenvalues of $D^2$, too.
The second eigenvalue can be estimated by
\bdm
a^2 n^2 \ = \ \mu_1(D^2) \ < \ \mu_2(D^2) \ \leq \ \Big(\sqrt{\lambda^0_1 \, + \, a^2 (1 - n)^2} \ - \ |a|\Big)^2 
\edm
Finally, if 
\bdm
a^2 n^2 \ = \ \mu_1(D^2)\ < \ \mu(D^2) \ < \  \Big(\sqrt{\lambda^0_1 \, + \, a^2 (1 - n)^2} \ - \ |a|\Big)^2 
\edm
is any ``small''  eigenvalue and $\psi^*$ the eigenspinor, then the inner product $\langle 
\psi \, , \, \psi^* \rangle$ vanishes identically.
\end{thm}
\vspace{5mm}

\noindent
Spinor fields $\psi^* = \eta \cdot \psi$ given as the Clifford product of the
Killing spinor $\psi$ by a $1$-form $\eta$ satisfy automatically the condition
$\langle \psi \, ,\, \psi^* \rangle = 0$. We compute the Dirac operator,
\bdm
D(\psi^*) \ = \ (n - 2) \, a \cdot \eta \cdot \psi \ + \ d \eta \cdot \psi \ + \
\delta\eta \cdot \psi \ .
\edm
A first application of the formula is the following
\begin{prop}
If the Killing spinor $\psi$ is preserved by the Killing $1$-form $\eta$,
$\mathcal{L}_{\eta}(\psi) = 0$, then the spinor $\psi^* = \eta \cdot \psi$ is an
eigenspinor,
\bdm
D(\psi^*) \ = \ (n \, + \, 2) \, a \cdot \psi^* \ .
\edm
\end{prop}
\begin{proof}
The $1$-form $\eta$ is coclosed, $\delta \eta \, = \, 0$, and the formula for
the Lie derivative
\bdm
0 \ = \ \mathcal{L}_{\eta}(\psi) \ = \ \nabla_{\eta}\psi \, - \, \frac{1}{4}
d\eta \cdot \psi \ = \ a \cdot \eta \cdot \psi \, - \,  \, \frac{1}{4}
d\eta \cdot \psi \ ,
\edm
see \cite{BourGau}, yields the result.
\end{proof}
\noindent
More generally, suppose that $\psi^*$ is an eigenspinor. The eigenvalue equation
$D(\psi^*) = m \cdot \psi^*$ reads as
\bdm
\Big\{ \big( (n-2) a \, - \, m \big) \, \eta \ + \ d \eta \Big\} \cdot \psi \
= \ 0 \ .
\edm
The latter equation implies that the $1$-form has to be coclosed and 
an eigenform of the
Hodge-Laplace operator $\Delta_1$. Indeed, we have
\begin{lem} \label{Lemma2}
Let $\eta$ be a $1$-form and $0 \neq c$ a constant such that
\bdm
\big( c \cdot \eta \ + \ d \eta \big)\cdot \psi \ = \ 0
\edm
holds. Then the $1$-form is a divergence-free eigenform of $\Delta_1$,
\bdm
\delta \eta \ = \ 0 \, , \quad \Delta_1(\eta) \ = \ c (c \ - \ (2n - 6) a)\, 
\eta \ .
\edm
Furthermore, the conditions $\delta \eta = 0$ and $d \eta \cdot \psi =
0$ imply $\Delta_1(\eta) = 0$ (the case of $c = 0$).
\end{lem}  
\begin{proof}
Fix an orthonormal frame $e_1,
\ldots , e_n$ on $M^n$. Differentiate the equation for $\eta$ again, 
use the Killing equation for $\psi$  and
contract via the Clifford multiplication. Then we obtain
\bdm
\sum_{i=1}^n \Big\{c \, e_i \cdot (\nabla_{e_i}\eta) \, + \, a \cdot c \,
e_i \cdot \eta \cdot e_i  \, + \, e_i
\cdot (\nabla_{e_i} d \eta)  \, + \, a \, e_i \cdot d \eta
\cdot e_i \Big\}\cdot \psi \ = \ 0 \ .
\edm
The algebraic relations in the Clifford algebra 
\bdm
\sum_{i=1}^n e_i \cdot \eta \cdot e_i \ = \ (n - 2) \, \eta \, , \quad
\sum_{i=1}^n e_i \cdot d \eta \cdot e_i \ = \ (4 - n)  \, d \eta \, ,
\edm
as well as the well
known
formulas
\bdm
\delta \xi \ = \ - \, \sum_{i=1}^n e_i \haken \nabla_{e_i} \xi \, , \quad
d \xi \ = \ \sum_{i=1}^n e_i \wedge \nabla_{e_i} \xi \, , \quad
\edm 
for any differential form $\xi$ yield
\bdm
\sum_{i=1}^n e_i \cdot \nabla_{e_i}\eta \ = \ \delta \eta \ + \ d \eta \, , \quad
\sum_{i=1}^n e_i \cdot \nabla_{e_i}d \eta  \ = \ \delta d \eta \ + d\, d \eta
\ .
\edm
Inserting the latter formulas we obtain
\bdm
\Big\{ c \, \delta \eta \ + \ (c \, + \, (4-n) a) \cdot d \eta \ + \ (n-2) a \cdot c \cdot \eta \ + \
\delta \, d \, \eta \Big\} \cdot \psi \ = \ 0 \ .
\edm
We multiply by the spinor $\psi$ and obtain
$c \, \delta \eta \, |\psi|^2 = 0$ ,
i.e. the $1$-form $\eta$ is coclosed, $\delta \eta = 0$.
Finally we  use again the equation we started with, $d \eta \cdot \psi = - c \,
\eta \cdot \psi$. Then 
\bdm
\big\{ c \, ( - \, c \, + \, (2n - 6) \, a) \cdot \eta \ + \ \delta \, d \, \eta \Big\}
\cdot \psi \ = \ 0 \ .
\edm
This is a Clifford product of a $1$-form by a spinor. Consequently,
the $1$-form has to be trivial and the result follows.
\end{proof}
\noindent
Let us introduce the eigenvalues $0 < \Lambda_1 < \Lambda_2 < \ldots $ as 
numbers such that
the problem
\bdm
\Delta_1(\eta) \ = \ \Lambda \, \eta \ , \quad \delta \eta \ = \ 0
\edm
has a non-trivial solution. In general, if $\eta$ is
$1$-form on an n-dimensional manifold, then
\bdm
||\nabla \eta||^2 \ \geq \ \frac{1}{2} || d \eta ||^2 \ + \ \frac{1}{n} ||
\delta \eta||^2
\edm 
holds (see \cite{GaMey}, page 270). 
Consider a coclosed eigenform  
$\Delta_1(\eta) = \Lambda \, \eta, \, \delta \eta = 0$ on an Einstein space. Then
the latter inequality as well as the Weitzenb\"ock formula for $1$-forms imply
the estimate 
\bdm
\Lambda_1 \ \geq \ \frac{2 R}{n} \  = \ 8 (n-1) a^2 \ .
\edm

\noindent
The existence of a non-trivial solution of the equation $(c \cdot \eta + d \eta) \cdot \psi =
0$ implies the inequality
\bdm
8 (n-1) a^2 \ \leq \ \Lambda_1 \ \leq \ c \big(c \, - \, (2n - 6) a \big) \ .
\edm
The latter inequality estimates the absolute value of $c$, 
\bdm
\sqrt{\Lambda_1 \, + \, a^2 (n - 3)^2} \ - \ (n -3) |a| \ \leq \ |c| \ .
\edm
Inserting $ m = (n - 2) a - c$, we obtain
\bdm
\sqrt{\Lambda_1 \, + \, a^2 (n - 3)^2} \ - \  |a| \ \leq \ |m| \ .
\edm
We summarize the result of the previous discussion.
\begin{thm}
The spinor field $\psi^* = \eta \cdot \psi$ is an eigenspinor, $D(\psi^*) =
m \, \psi^*$, if and only if
\bdm
\Big\{ \big( (n-2) a \, - \, m \big) \, \eta \ + \ d \eta \Big\} \cdot \psi \
= \ 0 \ .
\edm
In this case the
$1$-form $\eta$ is a coclosed eigenform of the Laplace operator, and the
eigenvalue
can be estimated by
\bdm
\sqrt{\Lambda_1 \, + \, a^2 (n - 3)^2} \ - \  |a| \ \leq \ |m| \ .
\edm
\end{thm}

\noindent
Since in dimension seven {\it any} real spinor field $\psi^*$ being orthogonal to
the Killing spinor $\psi$ is defined by a unique $1$-form $\eta$, we are able
to estimate the second eigenvalue of $D^2$ on $7$-dimensional manifolds $M^7$ with
Killing spinors. However, it may happen that $M^7$ admits more than only one 
Killing spinor. Consequently, we have to investigate the  equation
$\big\{ \big( (n-2) a \, - \, m \big) \, \eta \ + \ d \eta \big\} \cdot \psi \
= \ 0 $ for the $1$-form $\eta$ in more detail.
\begin{exa}
Let $M^5 \neq S^5$ be a Einstein-Sasaki manifold of dimension five. Then $a = 
 \pm 1/2$ and
$M^5$ admits exactly one $(\pm 1/2)$-Killing spinor. We obtain
\bdm
\frac{25}{4} \ = \ \mu_1(D^2) \ < \ \mu_2(D^2) \ \leq \ \Big(\sqrt{\lambda^0_1
  \, + \, 4} \ - \ \frac{1}{2}\Big)^2 \ .
\edm
If the dimension of the isometry group
is at least two, then there exists a Killing vector field $X$ preserving
the Killing spinor.
The spinor field $\psi^* := X \cdot \psi$ solves the Dirac equation with
eigenvalue $m = 7/2$ and we obtain the upper bound
\bdm
\mu_2(D^2) \ \leq \ \frac{49}{4} \ .
\edm
The first Laplace eigenvalue $\lambda^0_1 \geq 5$ has been computed for
special families, see \cite{Gibbons} and \cite{KSY}. There are examples with a
$2$-dimensional isometry group and
\bdm
\lambda^0_1 \ = \ \frac{33}{4} \, , \quad \mu_2(D^2) \ \leq \ 9 \ .
\edm
\end{exa}
\section{The Dirac spectrum of a nearly parallel $G_2$-manifold}\noindent
%

\noindent
A $7$-dimensional, simply-connected Riemannian spin manifold admits a Killing
spinor if and only if there is a (nearly)-parallel $G_2$-structure $\omega^3$
in the sense of \cite{FerGray}, 
\bdm
d \, \omega^3 \ = \ - 8 \, a \, * \omega^3 \, , \quad \delta \omega^3 \ = \ 0
\ .
\edm 
The $3$-form $\omega^3$ is defined by the spinor in a unique way, see
\cite{FKMS}. If $a=0$, then $M^7$ has a parallel $G_2$-structure and compact examples are
known, see \cite{Joyce}. If $a \neq 0$ and $M^7 \neq S^7$, then there are
three types. Indeed, denote by $m_a$ the dimension of the space of
all Killing spinors. Then $1 \leq m_a \leq 3$, $m_{-a} = 0$  and $M^7$ is either a 
$3$-Sasakian manifold ($m_a = 3$), a Sasaki-Einstein manifold ($m_a = 2)$ or
a proper nearly parallel manifold ($m_a = 1$), see \cite{FrKa1},
\cite{FrKa2}. 
Compact examples for any type are known, see \cite{BG}, \cite{FKMS}.\\

\noindent
The first
positive eigenvalue of the Laplace operator is bounded by ($a \neq 0$)
\bdm
\frac{R}{6} \ = \ 28 \, a^2 \ \leq \ \lambda^0_1 \ < \ \lambda^0_2 \ < \ldots
\edm
 and equality occurs if and only of $M^7$ is isometric to the sphere $S^7$
 (Lichnerowicz-Obata theorem). The further invariants we need are the numbers
\bdm
48 a^2 \ \leq \ \Lambda_1 \ \leq  \ \lambda^1_{1,\pm} \ <   \lambda^1_{2,\pm} \ < \ldots 
\edm
such that there exists a $1$-form $\eta$ with
\bdm
\Delta_1( \eta) \ = \  \lambda \, \eta \, , \quad    \delta \eta \ =
\ 0  \, , \quad \big(4 \, a \ \pm \ 
\sqrt{ 16 a^2 \, + \,  \lambda} \big)
\, \eta \ = \ - * \big( d \eta \wedge * \omega^3 \big) \ . 
\edm
This set is contained in the spectrum of the Laplace operator on
$1$-forms. Moreover, we have $\Lambda_i \leq \lambda^1_{i, \pm}$.
\begin{NB}
The equation $\big(4 \, a \ \pm \ 
\sqrt{ 16 a^2 \, + \,  \lambda} \big)
\, \eta \ = \ - * \big( d \eta \wedge * \omega^3 \big) $ implies $\delta \eta
= 0$ and $\Delta_1(\eta) = \lambda \, \eta$. Indeed, since $d * \omega^3 = 0$
we obtain
\bdm
\big(4 \, a \ \pm \ 
\sqrt{ 16 a^2 \, + \,  \lambda} \big)
\, d * \eta \ = \ - d * * \big( d \eta \wedge * \omega^3 \big) \ = \ 0 \ .
\edm
\end{NB}

\noindent
If $a=0$ and $ \lambda^1_{1,\pm} = 0$, then the Ricci flat manifold $M^7$
admits a parallel vector field, in particular $b_1(M^7) > 0$ holds. Consequently, for
simply-connected, compact and parallel $G_2$-manifolds, the numbers
$ \lambda^1_{i,\pm} > 0$ are positive.\\

\noindent
We formulate the main result of the section.
\begin{thm} \label{theoremn=7}
Let $(M^7, g)$ be a compact Riemannian spin manifold with a Killing
spinor. Then the spectrum of the Dirac operator consist of $(- 7\, a)$ and the following
sequences:
\bdm
- \, a \, \pm \, \sqrt{ 36 \, a^2 \ + \ \lambda^0_i } \, , \quad \lambda^0_i \,
\in \, \mathrm{Spec}(\Delta_0) \ .
\edm
\bdm
a \, - \, \sqrt{ 16 \, a^2 \ + \ \lambda^1_{i, +} } \, , \quad 
\mathrm{and}  \quad
a \, + \, \sqrt{ 16 \, a^2 \ + \ \lambda^1_{i, -} } \, , \quad 
i \ = \ 1 \, , \, 2 \, , \ldots \ .
\edm
\end{thm}
\vspace{5mm}

\noindent
The formulas simplify in case of a parallel $G_2$-structure ($a=0$).
\begin{proof}
The $7$-dimensional spin representation is real. Therefore, we consider real
spinor fields. The Killing spinor $\psi$ has constant length one, 
\bdm
\nabla_X \psi = \ a \cdot X \cdot \psi \, , \quad D(\psi) \ = \ - \, 7 a \cdot
\psi \, , \quad a \in \R^1 \ .
\edm 
It defines in a unique way a generic $3$-form $\omega^3$ and the Killing
equation as well as the link between $\psi$ and $\omega^3$ reads as
\bdm
d \omega^3 \ = \ - \, 8 a \, *\omega^3 \, , \quad \delta \omega^3 \ = \ 0 \, ,
\quad \omega^3 \cdot \psi \ = \ - \, 7 \cdot \psi \ ,
\edm
for details see \cite{FKMS}. $M^7$ is an Einstein space,
\bdm
\mathrm{Ric} \ = \ 24 \cdot a^2 \, \mathrm{Id} \, , \quad R \ = \ 7 \cdot 24 \cdot
a^2 \ .
\edm 
A purely algebraic computation in the $7$-dimensional spin representation
yields the following result.
\begin{lem} \label{Lemma1}
Let $\eta$ be a $1$-form, $\sigma$ a $2$-form and $c \in \R^1$ a real number.
Then
\bdm
\big( \eta \ + \ \sigma \ + \ c )\cdot \psi \ = \ 0
\edm 
is equivalent to
\bdm
\eta \ = \ - \, *  \big( \sigma \wedge * \omega^3 \big) \, \quad \mathrm{and}
\quad c \ = \ 0 \ .
\edm
\end{lem}

\noindent
Any real spinor field $\psi^*$ is given by a pair $(f , \eta)$  of a
real-valued function $f$ and a $1$-form $\eta$,
\bdm
\psi^* \ = \ f \cdot \psi \ + \ \eta \cdot \psi \ .
\edm
The eigenvalue
equation $D(\psi^*) = m \, \psi^*$ is equivalent to
\bdm
\Big( (- 7  a f \, - \, m f \, + \, \delta \eta) \ + \ 
(df \, + \, (5 a \, - \, m) \, \eta )\ + \ d \eta \Big) \cdot \psi \ = \ 0 \ .
\edm
By the algebraic Lemma \ref{Lemma1}, the Dirac equation reads now as
\bdm
(7  a  \, + \, m) f \ = \ \delta \eta \, , \quad 
df \, + \, (5a \, - \, m) \eta \ = \ - \, *\big(d\eta \wedge *\omega^3
\big) \ . 
\edm
Moreover, we know that $ m^2 \geq 49 a^2$ holds, see \cite{Fri1}. We
differentiate this system of first order partial differential equations
again and
we obtain the necessary condition
\bdm
\Delta_0(f) \ = \ (m \, - \, 5a)( 7a \, + \, m) f \ ,
\edm
i.e. the function $f$ is an eigenfunction of the Laplace operator.\\

\noindent
Let us first discuss the case $m = - \, 7a$. Then $f$ is constant and the
eigenspace $E_{-7a}(D)$ of the Dirac operator coincides with the space of all
Killing spinors with Killing number $a$, see \cite{Fri1}. This space becomes
isomorphic to
\bdm
E_{-7a}(D) \ = \ \R^1 \cdot \psi \ \oplus \ \big\{ \eta \, : \, 
12 a \, \eta \ = \ - \, *(d \eta \wedge *\omega^3) \quad \mathrm{and} \quad
\delta \eta \ = \ 0 \big\} \ .
\edm  
One needs the second equation only if $a = 0$.\\

\noindent
Suppose from now on that $m \neq - \, 7a$. Since $m^2 \geq 49 a^2$ we have
$(m - 5a)(7a + m) \neq 0$, too. The function $f$ is either non-trivial or
$f \equiv 0$. If $f \not\equiv 0$ then
\bdm
(m - 5a)(7a + m) \ = \ \lambda^0_i \ \in \ \mathrm{Spec}(\Delta_0) 
\edm
is a positive eigenvalue of the Laplace operator and
\bdm
m \ = \ - \, a \ \pm \, \sqrt{ 36 a^2 \, + \, \lambda^0_i} \ .
\edm
Conversely, given a non-trivial function with $\Delta_0(f) = \lambda^0_i f$ and 
$(m - 5a)(7a + m) = \lambda^0_i$, then the pair
\bdm
(f \, , \, \eta) \ := \ (f \, , \, \frac{1}{m \, - \ 5a} df )
\edm
defines a solution of the Dirac equation $D(\psi^*) = m \, \psi^*$. Any other
solution with the same function $f$ is given by a $1$-form $\eta_1$ being a
solution of the system
\bdm
\delta \eta_1 \ = \ 0 \, , \quad (5a \, - \, m) \eta_1 \ = \ - \, *(d \eta_1
\wedge * \omega^3) \ .
\edm
The latter equations describe the solutions of the Dirac equation for $f
\equiv 0$. \\

\noindent
Let us summarize the result. The eigenspace $E_{m}(D)$ consist
of all pairs
\bdm
( f \, , \, \frac{1}{m - 5a} df \ + \ \eta)
\edm
where $f \equiv 0$ or $f$ is an eigenfunction of the Laplace operator
\bdm
\Delta_0(f) \ = \ (m \, - \, 5a)(7a \, + \, m) f \ = \ \lambda^0_i f \, ,
\quad m \ = \ - \, a \, \pm \sqrt{36 a^2 \, + \, \lambda^0_i} \ . 
\edm
and $\eta$ is a special $1$-eigenform, 
\begin{eqnarray*}
(5a \, - \, m) \eta &=& - \, *(d \eta \wedge *\omega^3) \, , \\
\Delta_1(\eta) &=& (3a \, + \, m)(m \, - \, 5a) \eta = \lambda^1_i \eta
\, ,   \quad m \ = \ a \, \pm \sqrt{16 a^2 \, + \, \lambda^1_{i, \mp}} \ .
\end{eqnarray*}
\end{proof}
\begin{cor}
Let $(M^7, g)$ be a compact Riemannian spin manifold with a parallel
spinor $\psi$ and denote by $\omega^3$ the associated parallel
$G_2$-structure. Then $\mu_1(D^2) = 0$ and the second eigenvalues of the square of
the Dirac operator is given by
\bdm
\mu_2(D^2) \ = \ \mathrm{min} \big( \lambda^0_1 \, , \, \lambda^1_{1,+} \, ,
\, \lambda^1_{1,-} \big) \ .
\edm
$\lambda^0_1 > 0$ is the first positive eigenvalue of the Laplace
operator on functions and $\lambda^1_{1, \pm}$ are the first positive numbers
such that 
\bdm
\pm \, \sqrt{ \lambda} \cdot \eta \ = \ - \, *(d \eta \wedge * \omega^3)
\edm
admits a non-trivial solution.
\end{cor}

\begin{cor}
If $M^7$ admits exactly one Killing spinor ($m_a = 1 , \, a >0$), then 
$\lambda^0_1 > 28 a^2$ , $\lambda^1_{1,+} > 48 a^2$ , and
\bdm
\mu_2(D^2) \ = \ \mathrm{min} \Big( \big( \sqrt{36 a^2 \, + \, \lambda^0_1} \,
- \, a \big)^2 \, , \, \big( \sqrt{16 a^2 \, + \, \lambda^1_{1,+}} \, - \,
a \big)^2 \, , \,  \big(\sqrt{16 a^2 \, + \, \lambda^1_{1,-}} \, + \,
a \big)^2 \Big)
\edm
Moreover, we have $(\sqrt{16 a^2 \, + \, \lambda^1_{1,-}} \, + \, a ) \, \geq \, 9 a$
and $ (-a \, - \, \sqrt{36 a^2 \, + \, \lambda^0_1} ) \, \leq - \, 9a$.
\end{cor}
\begin{cor}
If $M^7$ admits exactly one Killing spinor ($m_a = 1 , \, a >0$) and at least
one non-trivial Killing vector field, then 
$\lambda^0_1 > 28 a^2$ , $\lambda^1_{1,+} > 48 a^2$ , $  \lambda^1_{1,-} = 48
a^2$. In particular,
\bdm
\mu_2(D^2) \ = \ \mathrm{min} \Big( \big( \sqrt{36 a^2 \, + \, \lambda^0_1} \,
- \, a \big)^2 \, , \, \big( \sqrt{16 a^2 \, + \, \lambda^1_{1,+}} \, - \,
a \big)^2 \, , \,  81 a^2) \ .
\edm
\end{cor}
\begin{proof}
Denote by $X$ the Killing vector field. Since the Killing spinor $\psi$ is
unique, we obtain
\bdm
0 \ = \ \mathcal{L}_X \psi \ = \ \nabla_X \psi \, - \, \frac{1}{4} \, d X \cdot
\psi \ = \ a \, X \cdot \psi \, - \, \frac{1}{4} \, d X \cdot
\psi \ ,
\edm
(see \cite{BourGau}). Then the spinor field $\psi^* := X \cdot \psi$ is an
eigenspinor, $D^2(\psi^*) = 81 a^2 \psi^*$.
\end{proof}

\begin{exa}
$\SO(5)/\SO_{ir}(3)$, $N(k,l) = \SU(3)/S^1_{k,l}$ and deformations
of $3$-Sasakian manifolds are examples with exactly one Killing spinor, see
\cite{FKMS}. 
\end{exa}

\begin{cor}
If $M^7 \neq S^7$ admits at least two Killing spinors, then $\lambda^0_1 > 28
a^2$,  $\lambda^1_{1,+} = 48 a^2$ and 
\bdm
\mu_2(D^2) \ = \ \mathrm{min} \Big( \big( \sqrt{36 a^2 \, + \, \lambda^0_1} \,
- \, a \big)^2 \, , \, \big( \sqrt{16 a^2 \, + \, \lambda^1_{2,+}} \, - \,
a \big)^2 \, , \,  \big(\sqrt{16 a^2 \, + \, \lambda^1_{1,-}} \, + \,
a \big)^2 \Big)
\edm 
\end{cor}

\begin{exa}
Two or three Killing spinors occur if $M^7$ is a
Sasaki-Einstein or a $3$-Sasaki manifold, see \cite{FrKa2}, \cite{FKMS}.
\end{exa}
%

\section{The $7$-dimensional  Sasaki-Einstein  case}\noindent
%
A simply-connected Sasaki-Einstein manifold $M^7 \neq S^7$ admits at least
 two Killing spinors with Killing number $a= 1/2$, see
\cite{FrKa2}. The scalar curvature equals $R = 42$ and and $\lambda^0_1 >
7$. The second Killing spinor $\psi^* = \eta \cdot \psi$ is given by a
$1$-form $\eta$ satisfying the equation
\bdm
6 \, \eta \ = \ - \, *(d \eta \wedge * \omega^3) \ .
\edm
and we obtain
\bdm
\lambda^1_{1,+} \ = \ 12 \, , \quad \lambda^1_{1,-} \ \geq \ 12 \ .
\edm
If $\lambda^1_{1, -} = 12$, then the corresponding eigenform
$\eta_1$ is a solution of the Laplace equation
\bdm
\Delta_1(\eta_1) \ = \ 12 \, \eta_1 \ = \ 2 \, \mathrm{Ric}(\eta_1) \, , \quad 
\delta \eta_1 \ = \ 0 \, ,
\edm
i.e. a Killing vector field (see \cite{Kobayashi}). Consequently, if the isometry group
of the Sasaki-Einstein manifold $M^7$ is one-dimensional, we have
$\lambda^1_{1,-} > 12$. If the dimension of the isometry group is at least
two, then there exists a Killing vector field $X$ preserving the Killing
spinor.
Then $\lambda^1_{1,-} = 12$ and and we obtain the following
\begin{thm}
Let $M^7 \neq S^7$ be a compact and simply-connected Sasaki-Einstein 
manifold and suppose that the dimension of the isometry group is at least
two. 
Then the second eigenvalue $\mu_2(D^2)$ is given by
\bdm
\frac{49}{4} \ = \  \mu_1(D^2) \ < \ \mu_2(D^2) \ = \ \mathrm{min} \Big( \Big( \sqrt{9 \, + \, \lambda^0_1} \,
- \, \frac{1}{2} \Big)^2 \, , \, \Big( \sqrt{4 \, + \, \lambda^1_{2,+}} \, - \,
 \frac{1}{2} \Big)^2 \, , \,  \frac{81}{4} \Big) \ .
\edm
\end{thm}
\vspace{3mm}

\begin{exa}
Let us discuss the case of $3$-Sasakian manifolds (see \cite{IK}, \cite{BG}). 
The isometry group is at
least $3$-dimensional. Moreover, there exists a spinor
field $\psi_0 = \eta_1 \cdot \psi$ such that
\bdm
D(\psi_0) \ = \ \frac{9}{2} \, \psi_0 \, , \quad |\psi_0| \ \equiv \, 1 
\edm
holds, see \cite{AgFr1} and \cite{Moro}. This spinor satisfies even a stronger
equation, namely
\bdm
\nabla_X \psi_0 \ = \ \frac{1}{2} X \cdot \psi_0 \quad \mathrm{if} \quad
X \in \mathrm{T}^v \, , \quad
\nabla_X \psi_0 \ = \ - \, \frac{3}{2} X \cdot \psi_0 \quad \mathrm{if} \quad
X \in \mathrm{T}^h \ .
\edm
The numbers
$\lambda^1_{1,+} = \lambda^1_{1,-}= 12$ coincide.
The solutions of the equations 
\bdm
 - \, 2 \, \eta \ = \ - \, *(d \eta \wedge * \omega^3) \, , \quad
\mbox{and} \quad  6 \, \eta \ = \ - \, *(d \eta \wedge * \omega^3)
\edm
can be seen directly. Indeed,
consider the three contact structures $\eta_1, \eta_2,\eta_3$  of the 
$3$-Sasakian manifold. Then
\begin{eqnarray*}
d \eta_1 &=& - \, 2\, ( \eta_{23} + \eta_{45} + \eta_{67})  , \\ 
d \eta_2 &=& \ \  \, 2\, ( \eta_{13} - \eta_{46} + \eta_{57})  , \\
d \eta_3 &=& - \, 2\, ( \eta_{12} + \eta_{47} + \eta_{56})  . 
\end{eqnarray*}
and
\begin{eqnarray*}
\omega^3 &:=& \frac{1}{2} \big( \eta_1 \wedge d \eta_1 \, - \,  
\eta_2 \wedge d \eta_2 \, - \,  \eta_3 \wedge d \eta_3\big) \\
*\omega^3 &=& - \, \frac{1}{8} \big( d \eta_1 \wedge d \eta_1 \, - \,  
d \eta_2 \wedge d \eta_2 \, - \,  d \eta_3 \wedge d \eta_3\big) 
\end{eqnarray*}
is one of the associated nearly parallel $G_2$-structures (see \cite{AgFr1}).
A purely algebraic
computation yields the relations
\bdm
- \, 2 \, \eta_1 \ = \ - \, *(d \eta_1 \wedge * \omega^3) \, , \quad
6 \, \eta_2 \ = \ - \, *(d \eta_2 \wedge * \omega^3) \, , \quad 
6 \, \eta_3 \ = \ - \, *(d \eta_3 \wedge * \omega^3) \ .
\edm
\end{exa}
\vspace{2mm}

\noindent
Let us consider the regular case. 
The contact structure induces a $\SO(2)$-action
and the orbit space $X^6 := M^7/\SO(2)$ is a 
$6$-dimensional K\"ahler-Einstein orbifold with scalar curvature $\bar{R} =
48$, see \cite{FrKa2}, \cite{BG}. The projection $\pi : M^7 \rightarrow
X^6$ is a Riemannian submersion with a totally geodesic fiber and
commutes with the Laplacian $\Delta_0$ on functions. 
Suppose that there exists an invariant eigenfunction.
It projects to an eigenfunction of
$\Delta_0$ on the orbifold $X^6$. Consequently, the corresponding 
eigenvalue (not the
multiplicity) of $M^7$ and $X^6$ coincides,
\bdm
\lambda^0_i(M^7) \ = \ \lambda^0_i(X^6) \ .
\edm  
In this case  we can apply an estimate proved by
A. Lichnerowicz for $\lambda^0_1$ of a smooth K\"ahler-Einstein manifold
(see \cite{Ballmann}, page 84)
\bdm
\lambda^0_1(X^6) \ \geq \ \frac{\bar{R}}{3} \ = \ 16 \ .
\edm 
Consequently, the general estimate of $\mu_2(D^2)$ simplifies.
\begin{thm}
Let $M^7 \neq S^7$ be a compact and simply-connected Sasaki-Einstein 
manifold. Suppose that the dimension of the isometry group is at least
two, that there exists a $\SO(2)$-invariant eigenfunction
with eigenvalue $\lambda^0_1$ and suppose that $X^6 = M^7/\SO(2)$ is smooth. 
Then the second eigenvalue $\mu_2(D^2)$ is given by
\bdm
\frac{49}{4} \ = \  \mu_1(D^2) \ < \ \mu_2(D^2) \ = \ \mathrm{min} \Big( \Big( \sqrt{4 \, + \, \lambda^1_{2,+}} \, - \,
 \frac{1}{2} \Big)^2 \, , \,  \frac{81}{4} \Big) \ .
\edm
\end{thm}
\begin{NB}
The assumptions are satisfied if $M^7/\SO(2)$ is smooth and there eixts a
subgroup $\SO(2) \subset \mathrm{G} \subset \mathrm{Iso}(M^7)$ being
isomorphic to $ G = \SO(3), \Spin(3)$. Indeed, any irreducible, real 
$\mathrm{G}$-representation admits a $\SO(2)$-invariant vector. 
\end{NB}

\vspace{1cm}

    
\end{document}